\documentclass[12pt]{article}

\begin{document}

\def \dU{{\cal U}(a,x)\,{\cal U}(a,-x)}
\title{
$$
{}
$$
$$
{}
$$
{\bf THE INTEGRAL REPRESENTATION}\\
{\bf FOR THE PRODUCT OF TWO}\\
{\bf PARABOLIC CYLINDER FUNCTIONS}
{\bf $D_\nu(x) D_\nu(-x)$ AT $Re\,\,\nu<0$ BY MEANS OF}\\
{\bf THE FUNDAMENTAL SOLUTION}\\
{\bf OF A LANDAU-TYPE OPERATOR }
}

\author{C.~Malyshev
\\
{\it V.A.Steklov Institute of Mathematics at St.-Petersburg}\\
{\it     Fontanka 27,~St.-Petersburg 191011, Russia}
\\ E-mail: malyshev@pdmi.ras.ru}

\maketitle

\def \text{\rm}
\def \hk{\widehat k}
\def \hl{\widehat l}
\def \wh{\widehat}
\def \wt{\widetilde}
\def \3{^3{\text He}-{\text A}}
\def \C{\text{C}}
\def \om{\omega}
\def \be{\beta}
\def \rot{{\text{rot}}}
\def \Dl{\Delta}
\def \dl{\delta}
\def \al{\alpha}
\def \ph{\phi}
\def \cd{\partial}
\def \si{\sigma}
\def \ep{\epsilon}
\def \ga{\gamma}
\def \Ga{\Gamma}
\def \Om{\Omega}
\def \la{\lambda}
\def \IM{\text{Im}}
\def \RE{\text{Re}}
\def \df{\frac{d}{\,dx\,}}
\def \dfz{\frac{d}{\,dz\,}}
\def \vj{\vec j}
\def \lrot{|\,\hl\times\rot\,\hl\,|}
\def \lll{\matrix{\, \cr ^{\stackrel{<}{\sim}}\cr}}
\def \rrr{\matrix{\, \cr ^{\stackrel{>}{\sim}}\cr}}
\def \BR{I\!\!R}
\def \BC{I\!\!\!\!C}
\def \UN{1\!\!{\text I}}
\def \x{{\text x}}
\def \d{{\text d}}

\begin{abstract}
The fundamental solution (Green's function) of a first order
matrix ordinary differential equation arising
in a Landau-type problem is calculated by two methods.
The coincidence of the two representations results in the integral
formula for the product of two parabolic cylinder functions
$D_\nu(x) D_\nu(-x)$ at ${\RE}\,\nu<0, x\in{\BR}$.
\end{abstract}

$$
{}
$$
$$
{}
$$
$$
{}
$$
\rightline{{\bf PDMI Preprint 04/2001}}
\rightline{{\bf MATH-CA/0106142}}

\newpage

\section{INTRODUCTION}

This paper is to point out another integral representation
for the product of two parabolic cylinder functions
$D_\nu(x) D_\nu(-x)$ at ${\RE}\,\nu<0, x\in{\BR}$. The
formula to be obtained is due to equivalence of two
representations for the fundamental solution (Green's
function) of a first order ($2\times 2$-matrix) ordinary
differential equation related to a couple of Weber's equations.
The parabolic cylinder functions have been extensively studied
in classical mathematical physics [1--5] (and refs. therein). In
particular, among various relations for them, the following
integral representations are known for products of two
parabolic cylinder functions:
$$
D_\nu(z)\,D_{-\nu-1}(z)\,, \eqno(1.1)
$$
($z$ real, ${\RE}\,\nu <0$, or ${\RE}\,z>0, {\RE}\,\nu >-1$),
$$
D_{-\nu-1} (z e^{i\frac\pi4})\,D_{-\nu-1}(z e^{-i\frac\pi4})\,,
\eqno(1.2)
$$
(${\RE}\,z>0, {\RE}\,\nu>-1$, or $|arg\,z|< \frac\pi4,
{\RE}\,\nu<-1$, or $-1 < {\RE}\,\nu < 0$),
$$
D_\nu(z e^{i\frac\pi4})\,D_\nu(z e^{-i\frac\pi4})\,,
\eqno(1.3)
$$
($|arg\,z|< \frac\pi4, {\RE}\,\nu<0$) [1, 2, 6]. Besides, the
following relations can be found in [5]:
$$
-\,\frac{\,1\,}z\,\dfz\left( D_\nu(z)\,D_\nu(-z)\right)\,=\,
\int_0^\infty J_0(z s) D^2_\nu(s) s d s\,,
\eqno(2)
$$
and [6]:
$$
(2z)^{-1}\left(D_\nu(-z)\,D_{\nu+1}(z)\,-
              \,D_{\nu+1}(-z)\,D_\nu(z)\right)
$$
$$
=\,(2z)^{-1}\left(z\,-\,\dfz \right) D_\nu(z)\,D_\nu(-z)
$$
$$
=\,\int_0^\infty J_0(z s) D_\nu(s)\,D_{\nu+1} (s) d s\,.
\eqno(3)
$$
Moreover, Ref.[1] suggests to consult with [7, 8, 9] for other
information concerning integral representations for the parabolic
cylinder functions. Essentially, apart from (1), only the
following combinations can be found in [7--16]:
$$
\left(D_\nu(z e^{i\frac\pi2})\,
               +\,D_\nu(z e^{-i\frac\pi2})\right)
\left(D_{-\nu-1}(z e^{i\frac\pi2})\,
               +\,D_{-\nu-1}(z e^{-i\frac\pi2})\right)
\eqno(4.1)
$$
($-1<{\RE}\,\nu< 0$),
$$
\left(D_\nu(z e^{i\frac\pi2})\,
               -\,D_\nu(z e^{-i\frac\pi2})\right)
\left(D_{-\nu-1}(z e^{i\frac\pi2})\,
               -\,D_{-\nu-1}(z e^{-i\frac\pi2})\right)
\eqno(4.2)
$$
($-2<{\RE}\,\nu< 0$),
$$
D_\nu(z)\,\left(D_\nu(z e^{i\frac\pi2})\,
              \pm\,D_\nu(z e^{-i\frac\pi2})\right)
\eqno(4.3)
$$
($|arg z|< \pi/4, {\RE}\,\nu< 0$)).
There are more other representations for products of $D_\nu(z)$
by means of integrals and series in [6, 17--24]. For instance,
the products at the same arguments can be found in [19, 21--23].
Certain representations through complex and contour integrals
are given in [22]. Expressions for products by means of
indefinite integrals are given in [19]. Recently, it has been
found that the wave field of the plane wave, which is scattered
by a cone of arbitrary shape, can also be expressed for singular
directions by means of the parabolic cylinder functions
(the scalar and electromagnetic cases) [25].

Combining (2) and (3) we obtain the following integral for the
product of two functions $D_\nu(z)$ at the opposite arguments:
$$
D_\nu(z)\,D_\nu(-z)\,=\,-2\,\int_0^\infty J_0(z s) D_\nu(s)
D^{\,\prime}_\nu(s) d s
\eqno(5)
$$
(notice that integration of (2) agrees with (5)). In its turn,
the present paper is to point out another formula for the
product $D_\nu(x)\,D_\nu(-x)$:
$$
D_\nu(x)\,D_\nu(-x)\,=\, 2^{-1/2} \Ga^{-1}(-\nu)
$$
$$
  \times\,\int_0^\infty \exp\left((\nu+\frac12) t-\frac{x^2}2\tanh
                     \frac t2\right)
\frac{dt}{{\sqrt{\sinh t}}}\,,
\qquad {\RE}\,\nu <0, x\in\BR\,, \eqno(6)
$$
which appears as a by-product of [26, 27] where a
theoretical problem of condensed matter physics is studied.
Specifically, the ground state average
of momentum operator (so-called, mass current) in weakly
inhomogeneous A-phase of superfluid helium-3 ($\3$)
is investigated in [26, 27]. Equation (6) appears since
the Green function of the quantum--statistical model under
consideration (i.e., the fundamental solution to an appropriate
matrix ordinary differential equation) can be obtained in two
equivalent forms: the integral [26], and the series [27].
Equating to each other, one obtains (6). The differential
equation under consideration below reminds the famous physical
problem of quantization for a spinning electron in a constant
homogeneous magnetic field (the Landau problem) [28].

It is especially convenient to apply to (6) the Laplace method
(steepest descent) in order to calculate asymptotical
expansions for the mass current [27]. Therefore, Eq.(6)
looks attractive since can be considered at the same
footing as (1)--(5), and, hopefully, can be useful for
applications, since contains only elementary functions under
integration. The present paper briefly recalls the main
points of [26] and [27] in Sections 2 and 3, respectively.
Section 4 contains special examples, and Section 5 concludes
the paper.

\section{THE INTEGRAL REPRESENTATION}

The master equation of our approach is given by the matrix
ordinary differential equation
$$
(i\om+{\cal H}(x))\,G(x)\,=\,e^{i x\xi}\,\UN\,,\qquad
{\cal H}(x)\,\equiv\,\left(\matrix{i\df & -(x+i\Dl)\cr
                       -(x-i\Dl) & -i\df\cr}\right),
\eqno(7)
$$
where $x\in\BR$ is a real variable;
$\Dl, \om, \xi$ are real parameters,
$i\equiv {\sqrt{ -1}}$, $\UN$ is unit matrix, and $G(x)$ is unknown
matrix. Equation (7) appears from the Dyson--Gorkov equation which
describes two-point correlation functions (Green's functions) of
$\3$ [26]. Under certain physically motivated conventions, Eq.(7)
appears from a spatially three-dimensional general Dyson--Gorkov
equation, and the variable $x$ turns out to be due to a separation
of the system into a collection of one-dimensional subsystems
(spherical coordinates are chosen). The parameter $\Dl$ depends,
generically, on the angle variables. The parameter $\om$ (7)
corresponds to a thermal Matsubara frequency. Solutions to the
Dyson--Gorkov equation describe the physics of the model in
question [26, 27].

It is convenient  to place the Dirac $\dl$-function in the R.H.S.
of (7):
$$
(i\om+{\cal H}(x))\,G(x,x^{\prime})\,=
                       \,e^{i x\xi}\dl(x-x^{\prime})\,\UN\,,
\eqno(8)
$$
so that (7) appears after integrating (8) over $x^{\prime}$ from
$-\infty$ to $\infty$ with $G(x)=\int G(x,x^{\prime}) dx^{\prime}$.
Let us put (8) in another equivalent form:
$$
(i\df + M (x))\,\wt G(x,x^{\prime})\,=
                       \,e^{i x\xi}\dl(x-x^{\prime})\,\UN\,,
\quad
M(x)\,\equiv\,\left(\matrix{i\om & x+i\Dl\cr
                       -x+i\Dl & -i\om\cr}\right).
\eqno(9)
$$
where $\wt G=\si G$, and $\si\equiv diag \{1, -1\}$ is the diagonal
matrix.

Following [26], let us assume that $G_0(x)$ is known which respects
$$
i \df G_0(x)\,=\,G_0(x) M(x).
\eqno(10)
$$
Being so, one is concerned, instead of (9), with solution of
$$\df(G_0(x) \wt G(x,x^{\prime}))\,=\,\frac 1i
                 e^{i x\xi}\dl(x-x^{\prime})\,G_0(x)\,.
\eqno(11)
$$
Without going into details of [26], let us calculate
$$
{\cal J}(x)=\int\!\!\!\int e^{-i x\xi}\,G(x,x^{\prime})
                                     dx^{\prime}d\xi
$$
($\xi$-integration is also from $-\infty$ to $\infty$). Notice
that $G(x,x^\prime)$ depends on $\xi$ due to (8).
Integrating (11) one obtains:
$$
{\cal J}_{1 1}(x)=\frac 1i(\det G_0)^{-1}
 \int \left[G_0^{22}(x)\int_{+\infty}^x\,G_0^{11}(s)\,e^{i\xi(s-x)}
 d s \right.
$$
$$
 \left.
-\,G_0^{12}(x)\int^x_{-\infty}\,G_0^{21}(s)e^{i\xi(s-x)}
 d s \right] d\xi\,,
\eqno(12)
$$
where the entries $G^{i j}_0$ respect (10). Equations (11), (12)
tell us that inhomogeneous Eq.(7) is solved, in fact, by the
method of variation of arbitrary constant [29].

In its turn, (10) can be reduced to the couple of second order
(Weber's, [3]) equations as follows. Let us note a useful matrix
relation:
$$
M^T\,-\,i\df\,=\,i\si\,u\, \left(\matrix{ a^+ &
\la \cr  \la^* & -a^- }\right) (\si\,u)^{-1},
\qquad u\,\equiv\,\frac1{{\sqrt 2}}\,
 \left(\matrix{ 1 & 1 \cr  i & -i }\right)\in U(2,{\BC})\,,
\eqno(13)
$$
where
$\si\,\equiv\,diag\{1, -1\}$, $a^{\pm}\equiv x\mp d/dx$,
$\la\equiv\om+i\Dl$, and superscripts $^T$ and $^*$ imply
transposed and complex conjugated, respectively. Taking $G_0$ in
the form
$$ G_0(x)\,=\,{\sqrt 2}
\left(\matrix{ h_1 & f_1 \cr h_2 & f_2 }\right) u^{-1}\,,
\eqno(14)
$$
and using
$$ \left(\matrix{ a^+ & \la \cr  \la^* & -a^- }\right) \,
\left(\matrix{ a^- & \la \cr  \la^* & -a^+ }\right)\,=\,
\left(\matrix{ a^+ a^-+|\la|^2 & 0 \cr 0 & a^- a^++|\la|^2}\right)
\,,
$$
one obtains unknown $f_{1, 2}\equiv f_{1, 2}(x)$ and
$h_{1, 2}\equiv h_{1, 2}(x)$ from the couple of Weber's
equations:
$$
\frac{d^2}{dy^2}\widetilde f
-\left(\frac{y^2}4-\frac12+|\lambda|^2\right)
\widetilde f=0\,,
\qquad
\frac{d^2}{dy^2}\widetilde h
-\left(\frac{y^2}4+\frac12+|\lambda|^2\right)
\widetilde h=0\,,
\eqno(14^{\prime})
$$
where $y\equiv x{\sqrt 2}$, the tilde implies that $f, h$ are
expressed through $y$, and 1, 2 are omitted. Besides, $\la$
is also rescaled in (14$^{\prime}$), i.e.,
$\la\equiv (\om+i\Dl)/{\sqrt 2}$ until the end of the paper.

Eventually, the following combinations for
$\widetilde h_{1,2}$ and $\widetilde f_{1,2}$
are chosen in [26] to express $G^{i j}_0$ which ensure the
integration in (12):
$$
 \widetilde h_1\,\pm\,\widetilde f_1\,=
\,{\cal U}_+ (y)\,\mp\,\la^{-1}{\cal U}_-(y)\,,
\quad \pm i (\widetilde f_2\,\mp\,\widetilde h_2)\,=
\,{\cal U}_+(-y)\,\mp\,\la^{-1}
{\cal U}_-(-y)\,,
$$
where ${\cal U}_{\pm}(y)\equiv {\cal U} (|\la^2|\pm 1/2, y)$ are the
parabolic cylinder functions (the notation ${\cal U}(a, y)$ implies
$D_{-a-\frac12}(y)$). The asymptotical behaviour of ${\cal U}(a,y)$
is given by the estimates [3]:
$$
{\cal U}(a,y)\simeq\exp(-y^2/4)\,y^{-a-1/2},\,\quad
{\cal U}(a,-y)\simeq\frac{(2\pi)^{1/2}}{
\Ga(a+1/2)}\,\exp(y^2/4)\,y^{a-1/2}\,,
$$
at $y\to +\infty$, $|y|\gg a$ (the leading terms).

The theory of distributions allows to calculate (12)
straightforwardly. For the given choice of $G^{i j}_0$, the order
of integrations in (12) can be changed at fixed $x$ (Fubini's
theorem). Integrating over $\xi$ first, we obtain $2\pi\dl (s-x)$.
Then, using the properties of $\dl$-function [29] we obtain:
$$
{\cal J}_{1 1}(x)\,=\,
i\pi \frac{ G_0^{11}(x)\,G_0^{22}(x)\,+\,
            G_0^{12}(x)\,G_0^{21}(x)}
          { G_0^{11}(x)\,G_0^{22}(x)\,-\,
            G_0^{12}(x)\,G_0^{21}(x)}
$$
(again in the $x$-notations).
A note about ${\cal J}_{2 1}$ will be given in the last
section. Other entries are not of interest for us.
The final answer for ${\cal J}_{1 1}(x)$ reads:
$$
{\cal J}_{1 1}(x)\,=\,\frac 1i \left(\frac\pi 2\right)^{1/2}
\left(\la\Ga(|\la|^2+1){\cal U}_+(x{\sqrt 2}) {\cal U}_+(-x{\sqrt 2})
\right.
$$
$$
\left. \,+\,
\la^* \Ga(|\la|^2){\cal U}_-(x{\sqrt 2}) {\cal U}_-(-x{\sqrt 2})
\right)\,,
\eqno(15)
$$
where $\Ga(\cdot)$ is the gamma function [2].

\section{THE SERIES REPRESENTATION}

Let us consider another approach to (7), (8) [27].
The operator ${\cal H}(x)$ (7) is self-adjoint on $\BR$, and it
can be conjugated by the constant matrix $u$ (13):
$$
u^{-1}\,{\cal H}\,u\,=\,{\cal H}_{em},\quad
{\cal H}_{em}\,=\,\left(\matrix{\Dl& i a^{-}\cr
-i a^+  &-\Dl\cr}\right).
$$
Squared operator ${\cal H}_{em}^2$ resembles the Hamiltonian
of a spinning electron in a constant homogeneous magnetic field
[28]. It is straightforward to obtain
eigenvalues $E_0$, $\pm E_n$ and eigenfunctions
$\widehat\Psi_0,\,\widehat\Psi^{\pm}_n\,(n\ge1)$
of ${\cal H}_{em}(x)$:
$$
\widehat\Psi_0\,=\,\left(\matrix {0\cr
\psi_0(x)\cr}\right),\quad E_0=-\Dl,
$$
$$
\widehat\Psi^{(s)}_n\,=\,\frac1{\sqrt {2 E_n}}
\left(\matrix {\,\,\,\,\,\sqrt{E_n+s\Dl\,}\,\,\psi_{n-1}(x) \cr
-is\sqrt{E_n-s\Dl\,}\,\,\psi_n(x)\cr}\right),\quad sE_n,
$$
where $s=\pm,\,\,E_n=\sqrt{\Dl^2+2n\,}$, and
$\psi_n(x)$ are the Hermite functions,
$\psi_n(x)\,=\,\pi^{-1/4}$
$\times (2^n n! )^{-1/2}\,e^{-x^2/2}\,H_n(x)$ [30].
Orthonormality and completness of the system $\widehat\Psi_0$,
$\widehat\Psi^{\pm}_n$ $(n\ge1)$ can be directly verified since
$\psi_n(x)$ $(n\ge 0)$ are orthogonal and complete in $L^2(\BR)$.

Let us take unknown $G$ in the form (14) and pass to the equation
$$
(i\om+{\cal H}_{em})\left(\matrix{ h\cr
f\cr}\right)\,=\,\dl(x-x')\,e^{ix\xi}
\left(\matrix{\frac{\,1\,}2\cr \frac{\,1\,}2\cr}\right)
$$
(the first column of the matrix equation in question).
Now the unknown $h$ and $f$ depend on $x$ and $x^{\prime}$,
and we expand $\left(\matrix{ h\cr f\cr}\right)$ in
$\widehat\Psi_0(x),\,\widehat\Psi^{\pm}_n(x)$
with the coefficients dependent on $x^{\prime}$:
$$
\left(\matrix{ h\cr f\cr}\right) (x, x^{\prime})\,=\,
b(x^{\prime})\,\widehat\Psi_0 (x)\,+\,
\sum_{s=+,-} \,\sum_{n=1}^{\infty} b^{(s)}_n(x^{\prime})
\,\widehat\Psi^{(s)}_n(x)\,.
$$
We define the coefficients [34] and, using
$$
G_{11}(x)=\int dx^{\prime} (h(x,x')\,+\,f(x,x'))\,,
$$
obtain ${\cal J}_{1 1}$ in the series form:
$$
{\cal J}_{1 1}(x)
\,=\,\frac1{i}\,\frac\pi{{\sqrt 2}}\,\sum\limits^\infty_{n=0}
\psi^2_n(x)\,\left(
\frac\la{|\la|^2+n+1}\,+\,\frac{\la^*}{|\la|^2+n}
\right).
\eqno(16)
$$

In [26, 27] one is concerned with (15) and (16)
themselves which correspond to the integrated entries of $G(x)$
(7). Let us now equate separately the real and
imaginary parts of (15) and (16). After the analytical continuation
$|\la|^2\longrightarrow z$, we deduce:
$$
\Ga(z)\,{\cal U}\left(z-\frac12, x{\sqrt 2}\right)
{\cal U}\left(z-\frac12,-x{\sqrt 2}\right)
\,=\,{\sqrt\pi}
\sum\limits^\infty_{n=0} \frac{\psi^2_n(x)}{n+z}\,,
\qquad {\RE}\,z>0\,.
\eqno(17)
$$
The series in R.H.S. of (17) is
convergent absolutely and uniformly on $\BR$, and it is
straightforward to express it in the integral form using
$$
(n+z)^{-1}\,=\,\int^\infty_0 e^{-t(n+z)} dt
\,,\qquad{\RE}\,z>0\,.
$$
Using the generating function of Hermite polynomials [2] to sum up
after exchanging summation and integration, one just obtains Eq.(6):
$$
\Ga(z)\,{\cal U}\left(z-\frac12, x\right) {\cal
U}\left(z-\frac12,-x\right)\,=
$$
$$
=\,\frac1{{\sqrt 2}} \int_0^\infty
            \exp\left(t(\frac12-z)-\frac{x^2}2\tanh \frac t2\right)
                          \frac{dt}{{\sqrt{\sinh t}}}\,,
\qquad {\RE}\,z>0, x\in\BR
\eqno(18.1)
$$
($x$ is rescaled, $x{\sqrt 2}\longrightarrow x$).
R.H.S. of (18.1) can be transformed to the integral
which is known [6, 31]:
$$
\int_0^1\frac{(1-s)^{z-1}}{(1+s)^z}\,\exp\left(-\frac{x^2}2 s
\right) \frac{d s}{{\sqrt s}}\,=\,
B\left(\frac12,z\right){\sf \Phi}_1\left(\frac12,z,z+\frac12;-1,
\frac{x^2}2\right)\,,\qquad {\RE}\,z>0\,,
\eqno(18.2)
$$
where ${\sf \Phi}_1(\al, \be, \ga; u, v)$ is the confluent
hypergeometric series in two variables [1, 32].

\section{THE SPECIAL CASES}

Two ways of solving (7) give us, naturally, two representations
for the same fundamental solution $G(x)$. Generally, it is a
powerful technique -- to equate equivalent representations of a
given Green function [33]. We shall consider three particular
examples in order to argue directly the interesting reductions of (18)
to the known facts. Equation (18.1) depends
on two parameters, $z\in{\BC}$ and $x\in{\BR}$, and their admissible
values constitute the domain ${\cal D}=$ $\{{\RE}\,z>0\}
\times{\BR}$. So, we shall consider the following three subsets in
${\cal D}$: the line, $z=1/2$; open half-plane, $x=0$;
asymptotic regions, $z\,\gg 1$, $|x|\ne 0$ is bounded, and
$|z|$ is limited, $|x|\gg 1$.

\subsection* {Example 1:\,$\quad z=\frac12, x\in {\BR}$.}

Using [2], we proceed:
$$
{\cal U}(0,x)\equiv D_{-\frac12}(x)=\left(\frac x{2\pi}\right)^{1/2}
              K_{\frac14}\left(\frac{x^2}4\right)\,,
$$
$$
{\cal U}(0,-x)\equiv D_{-\frac12}(-x)=-i\left(
\frac x{2\pi}\right)^{1/2} K_{\frac14}\left(\frac{x^2}4\right)
\,+\,
i \left(\frac x{\pi}\right)^{1/2}
                   K_{\frac14}\left(-\frac{x^2}4\right)\,.
$$
In their turn,
$$
K_{\frac14}\left(\frac{x^2}4\right)\,=\,\frac\pi{{\sqrt 2}}
\left(
I_{-\frac14}\left(\frac{x^2}4\right)
\,-\,I_{\frac14}\left(\frac{x^2}4\right)\right)\,,
$$

$$
K_{\frac14}\left(-\frac{x^2}4\right)\,=\,\frac\pi 2
\left(
(1-i)\,I_{-\frac14}\left(\frac{x^2}4\right)
\,-\,(1+i)\,I_{\frac14}\left(\frac{x^2}4\right)\right)\,.
$$
Then,
$$
{\cal U}(0,-x)\,=\,\frac{(\pi x)^{1/2}}2
\left(
I_{-\frac14}\left(\frac{x^2}4\right)
\,+\,I_{\frac14}\left(\frac{x^2}4\right)\right)\,,
$$
and
$$
\pi^{1/2}{\cal U}(0,x){\cal U}(0,-x)
\,=\,\frac{{\pi}^{1/2} x}{2^{3/2}}\,
           K_{\frac14}\left(\frac{x^2}4\right)
\left(I_{\frac14}\left(\frac{x^2}4\right)
      \,+\,I_{-\frac14}\left(\frac{x^2}4\right)\right)\,.
\eqno(19)
$$
On the other hand, the L.H.S. of (18.2) can be written at
$z=\frac12$ as
$$
\int^1_0 \exp\left(-\frac{x^2}2 s\right)\,
\frac{ds}{{\sqrt {s(1-s^2)}}}\,.
\eqno(20)
$$
Accordingly to [31] (Chapter 2), the integral
(20) is just equal to (19) (see APPENDIX).

\subsection* {Example 2:\,$\quad {\RE}\,z >0, x=0$.}

At $x=0$ the L.H.S. of (18.1) is [2]:
$$
\Ga (z) {\cal U}^2 (z-\frac12,0)\,\equiv\,
\Ga (z) D^2_{-z} (0)\,=\,
\frac{{\sqrt \pi}}2
\frac{\Ga \left(\frac z2\right)}
               {\Ga \left(\frac{z+1} 2\right)}\,.
$$
R.H.S. of (18.1) takes at $x=0$ the form of the definition
of the beta function [2]:
$$
\frac12\,\int^1_0 s^{z/2-1}(1-s)^{-1/2} d s \,=\,
\frac12\,B\left(\frac z2,\frac12\right)\,,
\qquad {\RE}\,\left(\frac z2\right)>0\,.
$$
The coincidence is clear.

\subsection* {Example 3:\,$\quad z \gg 1, |x|\ne 0$ is bounded,
                            or $|z|$ is bounded, $|x|\gg 1$.}

In leading order, we estimate the triple product
$\Ga(z){\cal U}(z-\frac12, x){\cal U}(z-\frac12,-x)$ in the L.H.S.
of (18.1) as [3]:
$$
\Ga(z)\,{\cal U}\left(z-\frac12, x\right) {\cal
U}\left(z-\frac12,-x\right)\,\simeq\, \left(\frac{2 \pi}{x^2+
4z-2}\right)^{1/2}\,, \eqno(21)
$$
where $z>1/2$, and $x^2+4z-2$ is
assumed positive and large.  By steepest descent we estimate the
R.H.S. of (18.1) at $z \gg 1$ ($|x|\ne 0$ is bounded) as $\left(\pi/(2
z)\right)^{1/2}$. Analogously, the integral in the L.H.S. of (18.2)
can be estimated at $|x|\gg 1$ ($|z|$ is bounded) as
$(2\pi)^{1/2}/|x|$.  In the both cases we get agreement with (21).
However, when $z$ and $|x|$ are large, (18.1) results in (21) as
well.

\def\tanhs{\text{\rm tanh}}
\def\sh{\text{\rm sh}}

\section{CONCLUSION}

To conclude, the following remarks are in order. The relation
$$
D_\nu(z)\,=\,2^{\frac\nu 2+\frac14}\frac1{{\sqrt z}}
W_{\frac\nu 2+\frac14,\pm\frac14} \left(\frac{z^2}2\right)
$$
is known which relates $D_\nu(z)$ to the Whittaker function.
Thus, our product $D_\nu(z) D_\nu(-z)$ can be expressed by
means of the integral representation for $W_{\kappa, \mu}(z)
W_{\chi, \mu}(z^{\prime})$ [2] (though containing the
hypergeometric function $_2F_1$ under integration).
However, some extra work is still required to relate (6) to the
corresponding formula. Consideration of the entry
${\cal J}_{2 1}$ does not give
really different formula: we obtain a representation which
can be deduced if one applies $x-d/dx$ to (6) (compare with
(3)). The representation (6) looks simpler in comparison with
analogous formulas mentioned above, and hopefully can be useful
in applications.

\subsection*{ACKNOWLEDGEMENTS}

The author is grateful to Prof. V. M. Babich for useful discussions.
The research described has been supported in part by RFBR,
No 01-01-01045

\subsection*{APPENDIX}

More precisely, [31] provides the formula for the integral
$$
\int_0^a \frac{x^{n-\frac12}}{{\sqrt{a^2-x^2}}}\,e^{-px}\,dx\,,
$$
which can be obtained from
$$
\int_0^a \frac{e^{-px}\,dx}{{\sqrt{x(a^2-x^2)}}}\,=
\,\frac{(\pi p)^{1/2}}{2}\,K_{\frac14}\left(\frac{a p}2\right)
$$
$$
\times\,
\left(I_{\frac14}\left(\frac{a p}2\right)
\,+\,I_{-\frac14}\left(\frac{a p}2\right)\right)\,,\quad
a>0\,,Re\,p>0
\eqno(a)
$$
by multiple differentiation $\partial^n/\partial (-p)^n$.
Besides, there exists another integral formula [31]:
$$
\int_a^\infty \frac{e^{-px}\,dx}{{\sqrt{x(x^2-a^2)}}}\,=
\,\frac{(\pi p)^{1/2}}{2}\,K_{\frac14}\left(\frac{a p}2\right)
$$
$$
\times\,
\left(I_{-\frac14}\left(\frac{a p}2\right)
\,-\,I_{\frac14}\left(\frac{a p}2\right)\right)\,,\quad
a>0\,,Re\,p>0\,.
\eqno(b)
$$
The couple of Eqs.(a) and (b) can be related to the
following couple of the integral transforms [35]:
$$
\int_0^\infty \frac{\sin{bx}\,dx}{{\sqrt{x(x^2+z^2)}}}\,=
\,\left(\frac{\pi b}2\right)^{1/2}\,K_{\frac14}\left(\frac{b z}2
\right)\,I_{\frac14}\left(\frac{b z}2\right)\,,\quad
b>0\,,Re\,z>0\,,
$$
$$
\int_0^\infty \frac{\cos{bx}\,dx}{{\sqrt{x(x^2+z^2)}}}\,=
\,\left(\frac{\pi b}2\right)^{1/2}\,K_{\frac14}\left(\frac{b z}2
\right) \,I_{-\frac14}\left(\frac{b z}2\right)\,,\quad
b>0\,,Re\,z>0\,.
$$

\subsection*{Notes added after release of the paper}

\subsubsection*{1.}
It is appropriate to point out that the representation
$$
\exp(-x^2)\,\left(H^2_\la (x) + G^2_\la (x)\right)\,=\,
2^{\la+1}\pi^{-1} \Ga(\la+1)
$$
$$
\times\int_0^\infty \exp\left(-(2\la +1)t+x^2\tanh t\right)
(\cosh t \sinh t)^{-1/2} dt
\eqno(N1)
$$
has been found in [36]. Here, $H_\la$ and $G_\la$ are the Hermite
functions which are related to the confluent hypergeometric
function $\Phi (a, b; z)$,
$$
\Phi (a, b; z)\,=\,\frac{\Ga (b)}{\Ga (a)}\,\sum_{n=0}^\infty
\frac{\Ga (a+n)}{\Ga (b+n)}\frac{z^n}{n!}\,,
$$
as follows [37]:
$$
2^{-\la} H_\la (x)\,=\,\frac{\Ga (1/2)}{\Ga (1/2-\la/2)}
\Phi (-\frac\la 2, \frac 12; x^2)\,+\,
                       \frac{\Ga (-1/2)}{\Ga (-\la/2)}\,x\,
\Phi (-\frac\la 2+\frac 12, \frac 32; x^2)\,,
$$
$$
2^{-\la} G_\la (x)\,=\,-\frac{\Ga (1/2)}{\Ga (1/2-\la/2)}
\tan\frac{\pi\la}2\,\Phi (-\frac\la 2, \frac 12; x^2)
+\,      \frac{\Ga (-1/2)}{\Ga (-\la/2)}\,x\,
\cot\frac{\pi\la}2\,\Phi (-\frac\la 2+\frac 12, \frac 32; x^2)\,.
$$
Equation (N1) appears in [36] as a confluent limit of an integral
representation for product of two Gegenbauer functions. The
representation (N1) turns out to be useful [36] in the problem of
monotonicity of differences of zeros of the Hermite functions [38,
39].

Using connection [37] between $\Phi (a, b; z)$ and the parabolic
cylinder functions, one finds:
$$
\exp(-x^2/2)\,H_\la (x)\,=\,2^{\la/2} D_\la (x{\sqrt 2})\,,
$$
$$
\eqno(N2)
$$
$$
\exp(-x^2/2)\,G_\la (x)\,=\,\frac{2^{\la/2}}{\sin\pi\la}
\left(\cos\pi\la D_\la (x{\sqrt 2})\,-\,D_\la (-x{\sqrt 2})
\right)\,.
$$
Eventually, using (N2) and the linear relations [2] between the
parabolic cylinder functions, it is straightforward to check that (N1)
and (6) (and thus (18)) are equivalent.

\subsubsection*{2.}
Equation analogous to (17) can also be written for the product of
two parabolic cylinder functions at coinciding arguments:
$$
\Ga(z)\,{\cal U}^2 \left(z-\frac12, x{\sqrt 2}\right)
\,=\,{\sqrt\pi}
\sum\limits^\infty_{n=0} \frac{(-1)^n\,\psi^2_n(x)}{n+z}\,,
\qquad {\RE}\,z>0\,,
\eqno(N3)
$$
where summation in R.H.S. goes over simple poles of $\Ga(z)$. In the
same way as in Sec.3, we obtain from (N3):
$$
D_{-z}^2(x)\,=\, 2^{-1/2} \Ga^{-1}(z)
\,\int_0^\infty \exp\left((\frac12-z) t-\frac{x^2}2\coth
                     \frac t2\right)
\frac{dt}{{\sqrt{\sinh t}}}\,,
\qquad {\RE}\,z >0\,,
\eqno(N4)
$$
or
$$
D_{-z}^2
(x)\,=\,\Ga^{-1}(z)\,\int_1^\infty\frac{(s-1)^{z-1}}{(s+1)^z}\,
\exp\left(-\frac{x^2}2 s\right) \frac{d s}{{\sqrt s}}\,.
\eqno(N5)
$$
However, the representation (N5) is known [31] just for
the square of the parabolic cylinder function but
in another form: i.e., for the integration variable obtained by
shift $s-1 \to s$. The integral in (18.2) and Eq.(N5) imply the
following combinations:
$$
D_{-z} (x)\left(D_{-z} (x) + D_{-z} (-x)\right)\,=\,
\Ga^{-1}(z)\,\int_0^\infty\frac{|s-1|^{z-1}}{(s+1)^z}\,
\exp\left(-\frac{x^2}2 s\right) \frac{d s}{{\sqrt s}}\,,
$$
$$
\eqno(N6)
$$
$$
D_{-z} (x)\left(D_{-z} (x) - D_{-z} (-x)\right)\,=\,
\Ga^{-1}(z)\,\int_0^\infty\frac{|s-1|^z}{(s+1)^z}\,\frac 1{s-1}
\exp\left(-\frac{x^2}2 s\right) \frac{d s}{{\sqrt s}}\,.
$$
Besides, (18.1) and (N4) can be re-written together in the following
form:
$$
D_{-z}(x)D_{-z}(\pm x)\,=
\,\frac{\,\exp(-x^2/2)\,}{2^{1/2} \Ga(z)}
\,\int_0^\infty \exp\left((\frac12-z) t\mp x^2 (e^t\mp 1)^{-1}
\right) \frac{dt}{{\sqrt{\sinh t}}}\,.
\eqno(N7)
$$

The Author is grateful to Prof. M. E. Muldoon for attracting his
attention to [36] and for interesting correspondence which inspired
to write down Eqs.(N3)--(N7) and to present the given notes.

\end{document}